\documentclass{article}
\usepackage[latin1]{inputenc}
\usepackage[english]{babel}
\usepackage{dsfont}
\usepackage[cyr]{aeguill}
\usepackage{amsmath,amssymb}
\usepackage{amsthm}
\usepackage{eurosym}
\usepackage{multicol}
\usepackage{color}
\usepackage{fancyhdr}
\usepackage[pdftex]{graphicx}
\usepackage[pdftex,colorlinks=automatic,linkcolor=black,citecolor=black,urlcolor=black]{hyperref}
\usepackage{verbatim}
\usepackage[linesnumbered,ruled,vlined]{algorithm2e}
\usepackage{stmaryrd}

\newcommand{\E}{\mathbb{E}} 

\theoremstyle{plain}
\newtheorem{thm}{Theorem}

\providecommand{\keywords}[1]{\textbf{\textit{Keywords --}} #1}

\SetKwInput{KwInput}{Input}                
\SetKwInput{KwOutput}{Output}              

\title{A comparison of maximum likelihood and absolute moments \\ for the estimation of Hurst exponents \\ in a stationary framework}

\author{Matthieu Garcin\footnote{ Léonard de Vinci Pôle Universitaire, Research center, 92916 Paris La Défense, France, matthieu.garcin@m4x.org.}}

\date{\today}

\begin{document}

\maketitle

\begin{abstract}
The absolute-moment method is widespread for estimating the Hurst exponent of a fractional Brownian motion $X$. But this method is biased when applied to a stationary version of $X$, in particular an inverse Lamperti transform of $X$, with a linear time contraction of parameter $\theta$. We present an adaptation of the absolute-moment method to this framework and we compare it to the maximum likelihood method, with simulations and an application to a financial time series. While it appears that the maximum-likelihood method is more accurate than the adapted absolute-moment estimation, this last method is not uninteresting for two reasons: it makes it possible to confirm visually that the model is well specified and it is computationally more performing.
\end{abstract}


\section{Introduction}

A fractional Brownian motion (fBm) $X$ is the only zero-mean Gaussian process with zero at the origin and with the following covariance function, for all $(s,t)\in\mathbb R^2$:
\begin{equation}\label{eq:CovFBM}
\E[X_tX_s]=\frac{\sigma^2}{2}(|t|^{2H}+|s|^{2H}-|t-s|^{2H}),
\end{equation}
where $H\in(0,1)$ and $\sigma>0$ are respectively the Hurst exponent of $X$ and its volatility parameter. The fBm was introduced by Mandelbrot and van Ness~\cite{MvN}. Beside the many extensions of this process~\cite{GarcinMulti}, several stationary specifications have been introduced~\cite{CKM,HN,Viitasaari}, with applications for example in meteorology~\cite{BSZ}, in fluid mechanics~\cite{Chevillard}, in finance~\cite{CV,GarcinLamperti}, or in medicine~\cite{SGKMBP}. For example, in finance, such adaptations of the fBm can be useful for modelling rates~\cite{Garcin2017,GarcinLamperti,GarcinMulti} or volatilities~\cite{GG,GJR}, which are reknown to be stationary. Among these stationary processes, we focus on the delampertized fBm~\cite{Lamperti,FBA,GarcinLamperti}. The Lamperti transform makes it possible to transform a self-similar process in a stationary process, as well as to do the reciprocal transformation. For example, the Ornstein-Uhlenbeck process, which is a widespread stationary process in quantitative finance, is the inverse Lamperti transform of a standard Brownian motion. We consider we are given a process $X$ which is $H$-self-similar, that is $\forall t,\lambda>0$, $X_t$ and $\lambda^{-H}X_{\lambda t}$ have the same probability distribution. For example an fBm of Hurst exponent $H$ is $H$-self-similar. According to the Lamperti transform, $X$ is derived from a stationary process $Y$. The process $X$ is the Lamperti transform of $Y$, $(\mathcal L_H Y)_t=t^{H}Y_{\ln(t)}$. The process $Y$ is the inverse Lamperti transform of $X$, also called delampertized $X$, $(\mathcal L^{-1}_H X)_t=\exp(-Ht)X_{\exp(t)}$. In what follows, we will consider the processes $Y_t=(\mathcal L^{-1}_{H,\theta} X)_{t}=(\mathcal L^{-1}_H X)_{\theta t}$ and $Z_t=(\mathcal L_{H',\theta'} Y)_{t}=t^{H'}Y_{\ln(t)/\theta'}$, in which we have added a parameter to linearly contract the time, which plays a role similar to the strength of a mean reversion~\cite{GarcinLamperti}.

We stress the fact that the delampertized fBm is not similar to the fractional Ornstein-Uhlenbeck (fOU) model, which is another widespread specification of stationary process using the fBm. The main difference is that, when $H>1/2$, the fOU model is a long-memory process, whereas the delampertized fBm is not~\cite{CKM,GarcinLamperti}.

Many methods enable to estimate the parameters of an fBm~\cite{MW}. The fOU process also has some specific estimation methods~\cite{BI,ST}. But only a few articles deal with the estimation of a delampertized fBm. We can cite, for example, an attempt using the covariance kernel of the process~\cite{MM}. As we are able to write the covariance of a delampertized fBm, we can estimate its parameters thanks to a maximum-likelihood (ML) method. However, like for other fractional processes~\cite{RL}, the ML estimation is time-consuming since it requires the inversion of the covariance matrix~\cite{MM}. For this reason, it may be useful to introduce other estimation methods, whose output may be used for initializing the ML method. 


A popular method to estimate the Hurst exponent of an fBm relies on the absolute moments of its increments and exploits its self-similarity property~\cite{PLV,Coeur2001,Coeur2005,Bianchi,Garcin2017}. A central vizualisation tool in this perspective is the log-log plot, that is the plot of the log absolute moments of the increments with respect to the log duration of the increment~\cite{GG,GJR}. Contrary to the ML approach, the absolute-moment method makes it possible to assess the good specification of the model. Indeed, provided that increments are shown to be stationary and Gaussian, only an fBm can lead to a linear log-log plot. But this promising method is shown to be biased in the stationary extension of the fBm~\cite{GarcinLamperti}. We thus propose an adapted absolute-moment (AAM) method to estimate both $H$ and $\theta$. The rationale of the AAM method consists in applying a proper direct Lamperti transform on the stationary data in order to use the standard absolute-moment method on the transformed time series. The main contribution of the present paper is Algorithm~\ref{algo:obj}, which implements the AAM method.

We compare the AAM method with the ML method. A simulation study will show that the ML method is in general slightly more accurate than the AAM method. However, the AAM method has two advantages over the ML. As exposed above, this method is based on an affine regression of logarithms of absolute moments of increments on logarithms of time scales, whose linear shape means that the model is well specified. On the contrary, the highest value of likelihood only means that the choice of parameters is optimal, not that the model is well specified. The second advantage of the AAM method is the fast computation. Finally, when working with stationary data, we prefer the ML method for accuracy, but we recommend to use first the AAM method to initialize properly the ML method in order to reduce its computational time. We can then apply the direct Lamperti transform with the ML-optimized parameters, following the AAM rationale, in order to observe whether the delampertized fBm model is well specified or not.

An application to the USD/EUR FX rate also shows the relevance of the model and of the estimation method in practice. Indeed, in finance, many studies have given to several kinds of time series the reputation to be stationary, like interest rates~\cite{CTY,ABKZ}, volatilities~\cite{CV,BLP,GG}, FX rates~\cite{NS,GarcinLamperti}, or even transaction volumes~\cite{Shi}.


In what follows, we present successively the ML method, the AAM approach, a simulation study, and a short financial application, all in a stationary framework.

\section{ML estimation}

In this section, we are considering a vector $\textbf{S}=(Y_{t_1},...,Y_{t_N})'$ of observations. The corresponding observation times are $\textbf{T}=(t_1,...,t_N)$. We propose to estimate the parameters $H$ and $\theta$ of a delampertized fBm by maximizing the following log-likelihood:
\begin{equation}\label{eq:LL}
L(\textbf{S};\theta,H)=\frac{1}{2}\ln\left(\det\left[\Sigma^{-1}\right]\right)-\frac{N}{2}\ln(2\pi)-\frac{1}{2}\textbf{S}'\Sigma^{-1}\textbf{S},
\end{equation}
where $\Sigma$ is the covariance matrix of a standard delampertized fBm~\cite{FBA}: 
\begin{equation}\label{eq:Sigma}
\Sigma_{ij}=\cosh\left(\theta H[t_j-t_i]\right)-2^{2H-1}\left|\sinh\left(\frac{\theta[t_j-t_i]}{2}\right)\right|^{2H}.
\end{equation}

As one can see in Figures~\ref{fig:ProblOptim3D} and~\ref{fig:ProblOptim2D}, the log-likelihood is a very smooth function of the parameters. The use of a heuristic optimization algorithm thus makes it possible to find quite rapidly values of $\theta$ and $H$ that are close to the optimum. For instance, in this work, we implemented the Nelder-Mead algorithm~\cite{NM}. However, we observe a ridge of high likelihoods in Figure~\ref{fig:ProblOptim3D}, stretched following mainly the $\theta$ axis. It may thus be difficult for the algorithm to choose among the possible pairs of parameters on this ridge, which all lead to very close likelihoods. Simulation results will confirm this difficulty, in particular for the estimation of $\theta$.

The likelihood above is the one of a standard delampertized fBm, that is with variance 1 and mean 0. We may also be interested in affine transformations of this standard process for practical applications, in which, for example in finance, $\textbf{S}$ is a vector of log-prices or a vector of log-volatility. Indeed, in this case, $\textbf{S}$ will more realistically be of the form $\textbf{S}=\mu \textbf{1} + \sigma \textbf{Y}$, where $\textbf{Y}$ is a vector of a standard delampertized fBm, $\mu\in\mathbb R$, and $\sigma>0$. A straightforward extension of the likelihood in this case is:
$$L'(\textbf{S};\theta,H,\mu,\sigma)=\frac{1}{2}\ln\left(\det\left[\Sigma^{-1}\right]\right)-\frac{N}{2}\ln\left(2\pi\sigma^2\right)-\frac{1}{2}\left(\textbf{S}'-\mu \textbf{1}'\right)\Sigma^{-1}\left(\textbf{S}-\mu \textbf{1}\right),$$
with the same $\Sigma$ as in equation~\eqref{eq:Sigma}. Such a framework would increase the dimension of the space of parameters in which we apply an optimization algorithm to maximize the log-likelihood. This would result in the elongation of the computational time in the ML estimation. We note that this would not be the case in the AAM approach, in which the estimation procedure of $H$ and $\theta$ is not sensitive to $\mu$ and $\sigma$. In what follows, however, we focus on the standard case $\mu=0$ and $\sigma=1$.



\section{Absolute-moment estimation}

In this section, we expose the absolute-moment estimation. In the case of second-order moments, it corresponds to the analysis of the variogram of the process. We present successively the basic method, where the process is an fBm, some specificities related to the stationarity of the delampertized fBm, an adaptation of this estimation method to this particular case of stationary process, and a pseudo-code detailing the estimation algorithm.

\subsection{Basic case}

If $X$ is an fBm of parameters $H$ and $\sigma^2$, it has the property of $H$-self-similarity: whatever $c>0$ and $\tau\geq 0$, $X_{\tau}\overset{d}{=}c^{-H}X_{c\tau}$, where $\overset{d}{=}$ means equality in finite-dimensional distributions. Since $X_0=0$, we also have $X_{\tau}-X_0\overset{d}{=}c^{-H}(X_{c\tau}-X_0)$. By stationarity of the increments of the fBm $X$, whatever $s>0$, we have $X_{s+\tau}-X_s\overset{d}{=}X_{\tau}-X_0\overset{d}{=}\tau^{H}X_{1}$. In particular, this increment has a variance equal to $\tau^{2H}\sigma^2$. We can use this property to estimate $H$. Indeed, the empirical variance of the increments of duration $\tau$, $M_{\tau}(X)$, is a sum of identically distributed terms, and it converges towards $\tau^{2H}\sigma^2$. The basic absolute-moment estimation thus consists in identifying the slope, expectedly equal to $2H$, of the log-log plot $\ln(\tau)\mapsto\ln(M_{\tau}(X))$~\cite{Coeur2005}.

We can replace the variance $M_{\tau}$ by an empirical absolute moment of another order~\cite{Coeur2005}. More generally, we can define an absolute moment of order $k>0$ for any process $S$ observed between times $t_a$ and $t_b$, with $N$ equispaced increments of duration $\tau=(t_b-t_a)/N$: 
$$M_{k,N,t_a,t_b}(S) = \frac{1}{N}\sum_{i=1}^{N}{\left|S_{t_a+(t_b-t_a)\times i/N}-S_{t_a+(t_b-t_a)\times(i-1)/N}\right|^k}.$$
The estimator of $H$ is then $-1/k$ times the slope of the log-log plot $\ln(N)\mapsto \ln(M_{k,N,t_a,t_b}(S))$.

\subsection{Stationarity and self-similarity}\label{sec:Statio}

We now observe a stationary process, $Y$. By Lamperti Theorem, it is the inverse Lamperti transform of parameters $H$ and $\theta$ of an $H$-self-similar process, $X$~\cite{FBA}. We also assume that $X$ is an fBm of Hurst exponent $H$. As the absolute-moment method is designed for the fBm, it is not the proper tool to estimate $H$ directly on $Y$~\cite{GarcinLamperti}. So, it sounds better to transform $Y$ first in a self-similar process $Z$, thanks to a direct Lamperti transform of parameters $H'$ and $\theta'$, and to estimate $H$ by applying the absolute-moment method on $Z$. We should preferably choose $H'=H$ and $\theta'=\theta$, so that $Z$ is equal to $X$ and is therefore an fBm, for which the absolute-moment method is relevant. But $H$ is unknown in practice and we cannot use it to transform $Y$ in $Z=X$. As a consequence, $Z$ is not necessarily an fBm, although it is $H'$-self-similar~\cite{FBA}. As exposed by the following theorem, increments of $Z$ are not stationary in general. It means that $M_{k,N,t_a,t_b}(Z)$ is not a sum of identically distributed terms, and that the log-log plot provides us with a slope which is not $-kH'$ in general. The log-log plot has even an asymptote of slope $-kH$ for $N\rightarrow +\infty$. These properties, exposed in Theorem~\ref{th:MomentIncr}, are pivotal for building the AAM estimation method.

\begin{thm}\label{th:MomentIncr}
Let $X$ be an fBm of Hurst exponent $H\in(0,1)$. Let $H'\in(0,1)$ and $\theta,\theta'>0$. Let $Y=\mathcal L_{H,\theta}^{-1}X$ and $Z=\mathcal L_{H',\theta'}Y$. Then:
\begin{itemize}
\item $X$ and $Y$ have stationary increments,
\item $Z$ has stationary increments if and only if $H'=H$ and $\theta'=\theta$.
\end{itemize}
In addition, $Z$ is such that, for $k,N\in\mathbb N$ and $t_b>t_a> 0$:
\begin{equation}\label{eq:th1NonAsympt}
\E\left[M_{k,N,t_a,t_b}(Z)\right]= \frac{A(\sigma,k)}{N}\sum_{i=1}^{N}\left[ t_{i+1}^{2H'} + t_i^{2H'} - \left(t_{i+1}t_i\right)^{H'-\frac{\theta}{\theta'}H}\left(t_{i+1}^{2H\theta/\theta'}+t_i^{2H\theta/\theta'}-\left[t_{i+1}^{\theta/\theta'}-t_i^{\theta/\theta'}\right]^{2H}\right) \right]^{k/2},
\end{equation}
where $A(\sigma,k)=\frac{2^{k/2}\Gamma\left(\frac{k+1}{2}\right)}{\Gamma\left(\frac{1}{2}\right)}\sigma^k$ and $t_i=t_a+(t_b-t_a)\times(i-1)/N$. Moreover, we have asymptotically, when $N\rightarrow+\infty$:
\begin{equation}\label{eq:th1Asympt}
\E\left[M_{k,N,t_a,t_b}(Z)\right] \overset{N\rightarrow+\infty}{\sim} A(\sigma,k) \frac{t_b^{k(H'-H)+1}-t_a^{k(H'-H)+1}}{k(H'-H)+1} \left(\frac{\theta}{\theta'}\right)^{kH}\left(t_b-t_a\right)^{kH-1}N^{-kH}.
\end{equation}
\end{thm}

The proof of Theorem~\ref{th:MomentIncr} is postponed in Appendix~\ref{sec:ProofMomentIncr}.

In addition, with the assumptions of Theorem~\ref{th:MomentIncr}, $Z$ is $H'$-self-similar, while $X$ is $H$-self-similar and $Y$ is not self-similar but stationary~\cite{FBA}. The particular case $(H',\theta')=(H,\theta)$ is of interest, as it is the only case for which the log-log plot of $Z$ is indeed affine of slope $-kH'$, as we will see in Theorem~\ref{th:LogPlotLin}. The AAM method presented in the next paragraph capitalizes on this property.

\subsection{The AAM method}\label{sec:AAM}

The absolute-moment method is widespread and efficient for estimating Hurst exponents in the case of self-similar processes. We want to adapt this method to the case of a stationary process. To do so, we transform the stationary process in a self-similar one thanks to the direct Lamperti transform. We can then apply the basic absolute-moment method to this transformed process. However, this simple idea raises two difficulties: one is about the choice of the parameters in the Lamperti transform, the other on the estimation of the moments. For simplicity, we now focus on second-order moments.

Regarding the first issue, in order to transform the stationary process $Y$ into a self-similar one, we need to choose a scaling parameter $H'$ and a time change parameter $\theta'$. We restrict the framework by assuming that the stationary process is itself the inverse Lamperti transform of an fBm with unknown parameters $H$ and $\theta$. We will show in Theorem~\ref{th:LogPlotLin} that one and only one choice for the pair $(H',\theta')$ can lead to an fBm of Hurst exponent $H'$: it is precisely $(H,\theta)$. As a consequence, we are going to select the pair $(H',\theta')$ so that $Z=\mathcal L_{H',\theta'}Y$ is revealed to be an fBm of Hurst exponent $H'$. The absolute-moment method makes it possible to check if we meet this property. Indeed, if the process is an fBm, the log-log plot must be affine, and if its Hurst exponent is $H'$, the slope of the log-log plot must be $2H'$. We thus propose an iterative optimization procedure, using Nelder-Mead algorithm as for the ML method. The purpose of this procedure is to find the pair $(H',\theta')$ minimizing the objective function $f_{\mathcal S}(H',\theta')=|\hat H_{H',\theta'}-H'|+|\hat \alpha_{H',\theta'}-1|$, where $\hat H_{H',\theta'}$ and $\hat \alpha_{H',\theta'}$ are respectively half the slope and a linearity parameter of the log-log plot of $\mathcal L_{H',\theta'}Y$ on a set of scales $\mathcal S$. More specifically, if $M_{H',\theta',\tau}$ is the moment\footnote{ If $(H',\theta')\neq(H,\theta)$, the increments of the process are not stationary according to Theorem~\ref{th:MomentIncr}. So the word moment is improper for these increments. For concision, we however keep this word to designate the average of moments of increments at various times.} of the increments of $\mathcal L_{H',\theta'}Y$ of duration $\tau$, we obtain $\hat H_{H',\theta'}$ and $\hat \alpha_{H',\theta'}$ by the following regression: $\ln\left(M_{H',\theta',\tau}\right) = \ln\left(M_{H',\theta',\tau_{min}}\right) + 2 \hat H_{H',\theta'} (\ln(\tau)-\ln(\tau_{min}))^{\hat \alpha_{H',\theta'}}$, where $\tau_{min}$ is the smallest time step observed in the log-log plot and $\tau\in\mathcal S$. Theoretically, this regression leads to $\hat H_{H',\theta'}=H'$ and $\hat \alpha_{H',\theta'}=1$ only for $(H',\theta')=(H,\theta)$, as stated in the following theorem, which justifies the relevance of the proposed estimation method.

\begin{thm}\label{th:LogPlotLin}
The pair $(H',\theta')$ reaches the theoretical minimum of the function $f_{\mathcal S}$ whatever $\mathcal S$, if and only if $(H',\theta')=(H,\theta)$.
\end{thm}

The proof of Theorem~\ref{th:LogPlotLin} is postponed in Appendix~\ref{sec:ProofLogPlotLin}.

The other challenge in the stationary adaptation of the absolute-moment method is about estimating properly moments of increments of the transformed process $\mathcal L_{H',\theta'}Y$. This transformation indeed distorts the times. We now observe the process at irregularly sampled dates. The ML method can nicely face this irregular sampling~\cite{HSL}. But the estimation of absolute moments is not as straightforward. In the current work, we have implemented a semi-parametric estimation of the absolute moments. We have used a kernel regression on the absolute increments of the observed process with a fractal correction of their scale. For estimating the absolute moment of increments of duration $\tau$, we have indeed taken into account various absolute increments of duration $d$ close to $\tau$, which we have multiplied by $(\tau/d)^{2H'}$.

\subsection{Algorithms}

We now introduce pseudo-codes which expose how the AAM method is to be implemented. We divide the estimation procedure in two steps. The first one is about optimizing the parameters $H$ and $\theta$ and the second one about the implementation of the objective function to be optimized. The first step is not very specific and is shared by the two approaches. Indeed, in the ML approach, instead of minimizing the objective function implemented in Algorithm~\ref{algo:obj}, one uses Algorithm~\ref{algo:opt} to maximize the log-likelihood detailed in equation~\eqref{eq:LL}.

\begin{algorithm}[htbp]\label{algo:opt}
\caption{Optimization engine for both estimators}
\DontPrintSemicolon
  
  \KwInput{$f$\tcp*{the objective function, either the log-likelihood or $f_{\mathcal S}$}}
  \KwData{observations $\textbf S$ and corresponding times $\textbf{T}$}
  
 \tcc{Initialization step}
$(H_1,\theta_1)=(0.45,25)$\;
$(H_2,\theta_2)=(0.55,28)$\;
$(H_3,\theta_3)=(0.50,35)$\;  
  
  \tcc{Iteration step}
  \While{$\sum_{i=2}^3{\left\{\left|1-\frac{H_i}{H_1}\right|+\left|1-\frac{\theta_i}{\theta_1}\right|\right\}}>0.001$}
   {
   		$(H_i,\theta_i)_{i\in\{1,2,3\}} = NM(f,\textbf{T},\textbf S,(H_i,\theta_i)_{i\in\{1,2,3\}})$ \;
   		\tcp*{$NM$ is the iteration step of Nelder-Mead algorithm}
   }  
\KwOutput{$(H_1,\theta_1)$}
\end{algorithm}

The optimization engine in Algorithm~\ref{algo:opt} is quite simple to use and can be applied to both estimation methods. Nevertheless, this algorithm may be improved. First, the parameters $H$ and $\theta$ are constrained to remain in the intervals $(0,1)$ and $(0,+\infty)$. In a basic version of Nelder-Mead, we can simply force the parameters to be in these intervals, by changing forbidden values by one close to the border, following Box's method~\cite{Box}. Alternative promising solutions for dealing with constraints include a random version of Box's method as well as a transform such as $x\in\mathbb R\mapsto\frac{1}{2}+\frac{1}{\pi}\arctan(x)\in(0,1)$ for $H$ and $x\in\mathbb R\mapsto\exp(x)\in(0,+\infty)$ for $\theta$~\cite{LF}.

In Algorithm~\ref{algo:opt}, the initialization of the simplex relies on some arbitrary values which are not too close in order not to stop prematurely the algorithm, and not too distant from each other in order not to cross the bounds after the first iteration. The values displayed for the initialization in Algorithm~\ref{algo:opt} are only indicative. Other reasonable deterministic or random initializations lead to similar results. We could also call three times this optimizer with three distinct initializations and then use these three solutions as the initialization of a last Nelder-Mead optimization.

The stopping criterion
$$\sum_{i=2}^3{\left|1-\frac{H_i}{H_1}\right|+\left|1-\frac{\theta_i}{\theta_1}\right|}\leq 0.001$$ 
is intended to depict the situation in which the vertices of the simplex are close enough to consider that a local optimum has been found. A criterion on the number of iterations may also be introduced.

The objective function in Algorithm~\ref{algo:opt} is quite simple in the ML case since it is related to the log-likelihood. The most time-consuming part for computing this likelihood is the inversion of the symmetric definite-positive matrix $\Sigma$. The objective function in the AAM approach is based on the analysis of the log-log plot of the Lamperti transform of the data. It is exposed in Algorithm~\ref{algo:obj}, which is the main methodological contribution of this paper.

\begin{algorithm}[htbp]\label{algo:obj}
\caption{Objective function of the AAM approach}
\DontPrintSemicolon
  
  \KwInput{$H'$ and $\theta'$, parameters to be tested, and $\textbf R$, vector of $n$ scales to be considered}
  \KwData{vector $\textbf S$ of $N$ observations and corresponding times $\textbf{T}$}
  
	\tcc{We operate a Lamperti transform of the observed process}  
	\For{$i\in\llbracket 1,N\rrbracket$}    
        { 
        	$\textbf{T'}[i]=\exp(\theta'\textbf{T}[i])$\;
        	$\textbf{S'}[i]=\exp(\theta'H'\textbf{T}[i])\textbf{S}[i]$
        }	
	
	\tcc{We determine the approximated moments of the increments}
	\For{$i\in\llbracket 1,n\rrbracket$}{
		$\tau=\textbf R[i]$\;
		$\textbf M[i]=0$\;
		\For{$j\in\llbracket 1,N-1\rrbracket$, $k\in\llbracket j,N\rrbracket$}    
        {
			$d=\textbf{T'}[k]-\textbf{T'}[j]$\;      	
        	$w=K(d-\tau)$\;
        	$W=W+w$\;
        	$\textbf M[i]=\textbf M[i]+w(\textbf{S'}[k]-\textbf{S'}[j])^2(\tau/d)^{2H'}$
        }	
        $\textbf M[i]=\textbf M[i]/W$
    }
	\tcc{We calculate $\hat H$ with a log-log regression}
	\For{$i\in\llbracket 1,n\rrbracket$}{
		$\textbf M[i]=\ln(\textbf M[i])$\;
		$\textbf R[i]=\ln(\textbf R[i])$\;
	}
	$\hat H=\frac{1}{2}Slope(\textbf M,\textbf R)$
	\tcp*{$Slope(Y,X)$ is the slope of the linear regression}
	\tcp*{of the $Y[i]$ against the $X[i]$, with intercept}
	
	\tcc{We determine $\hat{\alpha}$ by regression}
	\For{$i\in\llbracket 1,n\rrbracket$}{
		$\textbf M[i]=\ln(\textbf M[i]-\textbf M[1])$\;
		$\textbf R[i]=\ln(\textbf R[i]-\textbf R[1])$\;
	}
	$\hat{\alpha}=Slope(\textbf M,\textbf R)$\;
	
  \KwOutput{$|1-\hat{\alpha}|+|H'-\hat H|$}
\end{algorithm}

We note that, in Algorithm~\ref{algo:obj}, we have simplified the two-dimensional regression leading to $\hat H_{H',\theta'}$ and $\hat \alpha_{H',\theta'}$ by two one-dimensional regressions. If $(H,\theta)=(H',\theta')$, the two regressions should lead to $\hat H=H'$ and $\hat{\alpha}=1$. This justifies the relevance of the simplification.

In the simulation study, we have considered equally-spaced observation times with a time step $h=t_{i+1}-t_i=0.001$. In the study presented in Section~\ref{sec:sim}, the set of scales $\mathcal S$, gathered in Algorithm~\ref{algo:obj} in the vector $\textbf R$, is such that the smallest scale considered is $\tau_{\min}=\textbf R[1]=\ln(e^{\theta h}-1)$ and the biggest $\textbf R[n]=\ln(e^{\theta h \rho N}-1)$, where $\rho\in(0,1)$ is a meta-parameter linking this biggest step to a proportion $\rho$ of the whole duration $N h$ of the signal. Other scales in $\mathcal S$ simply make a regular sampling of the interval $[\tau_{\min},\textbf R[n]]$.

The kernel $K$ used in the semi-parametric approximation of the moments does not play a decisive role in the results, according to the simulations we made. For decreasing the computational time, it is better to consider a kernel with bounded support, such as Epanechnikov or a truncated Gaussian. In this case, the choice of the bandwidth parameter of the kernel is to be made so that the total weight $W$, appearing in Algorithm~\ref{algo:obj} and corresponding to the scales selected in $\mathcal S$, is never equal to zero.

\section{Simulation study}\label{sec:sim}

We compute simulations and compare the ability of the ML and the AAM methods to estimate $H$ and $\theta$. We first illustrate both approaches in a simple case. We then analyse the impact of the parameters $H$ and $\theta$, as well as the number of observations $N$, on the performance of each method. Finally, we illustrate the benefit of using the AAM framework to determine whether the model is well specified or not.

\subsection{Illustration of the methods}

\begin{figure}[htbp]
	\centering
		\includegraphics[width=0.45\textwidth]{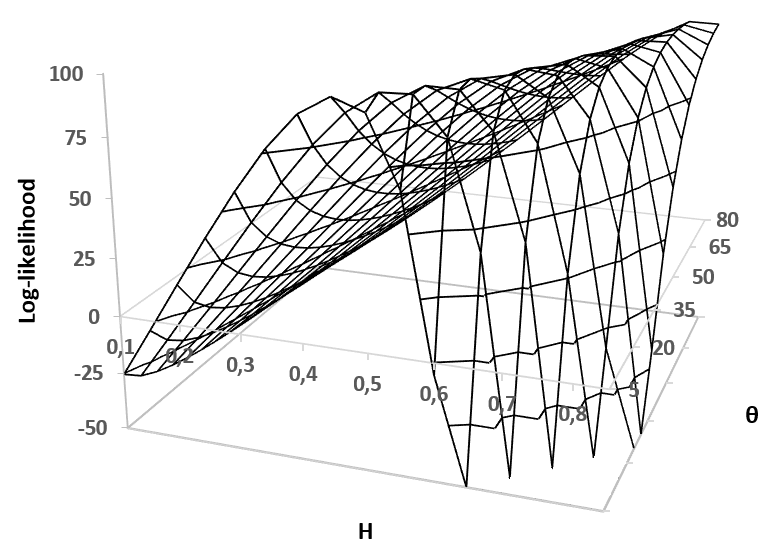}
		\includegraphics[width=0.45\textwidth]{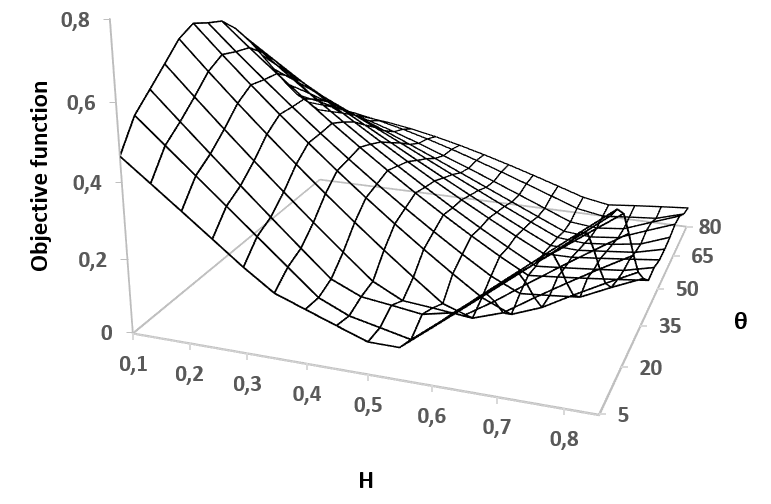}
\begin{minipage}{0.7\textwidth}\caption{Log-likelihood (left) and objective function $f_{\mathcal S}(H,\theta)$  (right) for various values of $H$ and $\theta$, on one simulation with the parameters 0.65 and 30.}
	\label{fig:ProblOptim3D}
\end{minipage}
\end{figure}

\begin{figure}[htbp]
	\centering
		\includegraphics[width=0.45\textwidth]{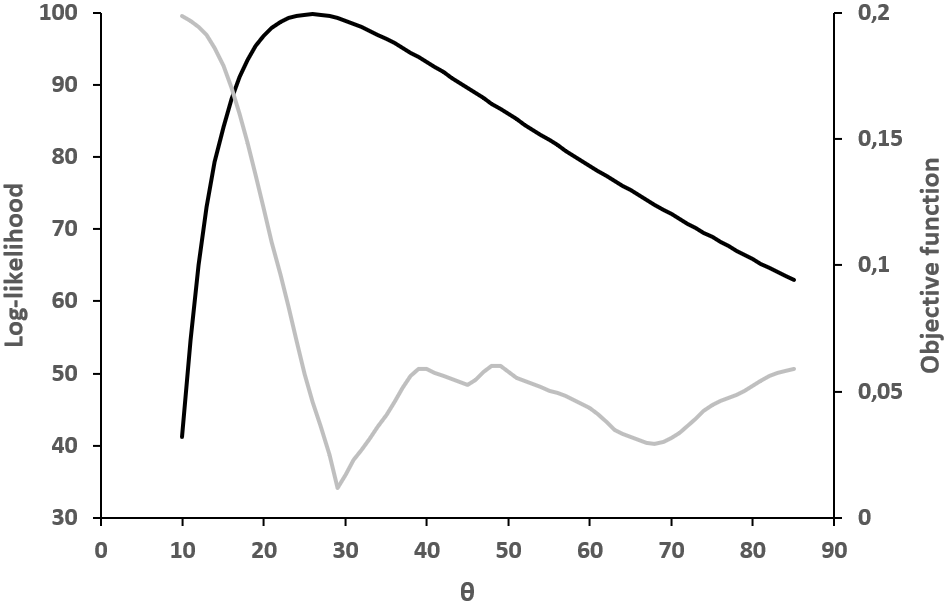}
		\includegraphics[width=0.45\textwidth]{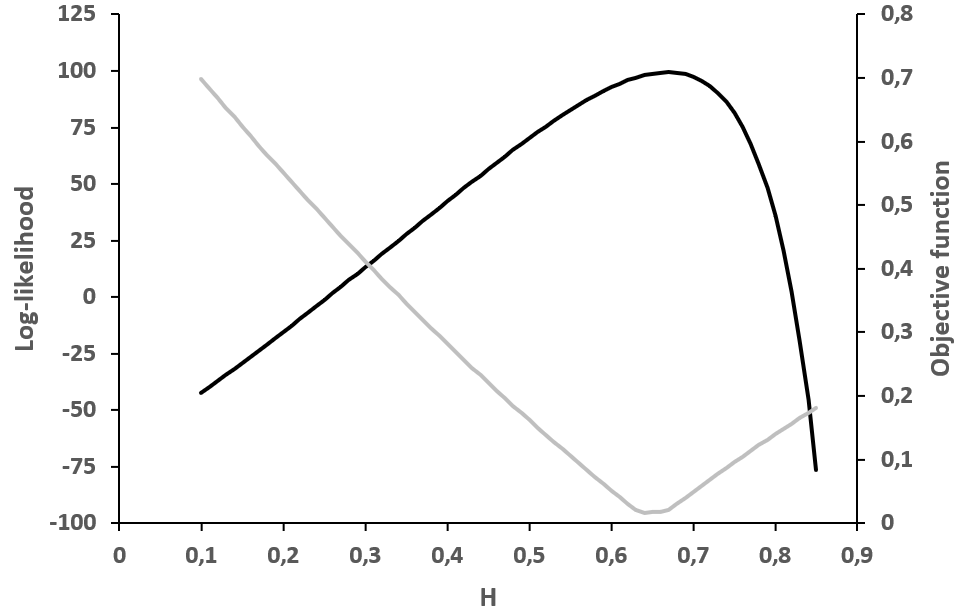}
\begin{minipage}{0.7\textwidth}\caption{Log-likelihood (black, left axis) and objective function $f_{\mathcal S}(H,\theta)$  (grey, right axis) for various values of $H$ and $\theta$, with fixed $H=0.65$ (left) or fixed $\theta=30$ (right), on one simulation with the parameters 0.65 and 30.}
	\label{fig:ProblOptim2D}
\end{minipage}
\end{figure}

We compare the ML and the AAM methods for estimating the parameters $H$ and $\theta$ of simulated delampertized fBms, with a Hurst exponent 0.65 and a time change parameter equal to 30. The functions that we optimize, either with a max for the likelihood or with a min for the objective function of the AAM method, are quite smooth curves, as one can see in Figure~\ref{fig:ProblOptim3D}. This eases the optimization. However, we see in Figure~\ref{fig:ProblOptim2D} that local minima may appear in the estimation of $\theta$ in the AAM approach. A more sophisticated algorithm, such as a genetic algorithm, instead of the Nelder-Mead algorithm, could probably increase the accuracy of the method. As already evoked, many pairs $(H,\theta)$ have a close likelihood, as a ridge with very high likelihoods appears in Figure~\ref{fig:ProblOptim3D}. This can make also the estimation in the ML approach unstable.

The average estimated parameters are close to the true parameters for both the ML and the AAM methods, as reported in Table~\ref{tab:simul}. The uncertainty in the estimation, revealed by the standard deviation, is higher for the ML with respect to $\theta$, and higher for the absolute moments with respect to $H$. 


\begin{table}[htb]
\centering
\begin{tabular}{|l|c|c|}
\hline
  Estimation method & Estimated $H$ & Estimated $\theta$ \\
\hline
ML & 0.659 \textit{(0.083)} & 36.7 \textit{(38.8)} \\
AAM & 0.651 \textit{(0.156)} & 34.6 \textit{(22.4)} \\
\hline
\end{tabular}
\begin{minipage}{0.7\textwidth}\caption{Average estimated parameters (and standard deviation) on 100 simulations with the parameters 0.65 and 30.}
\label{tab:simul}
\end{minipage}
\end{table}

\subsection{Impact of various parameters on the estimation}

We now study by simulations the impact of the input parameters and of methodological choices.

First, regarding the estimation of $\theta$ the ML method is undoubtedly more accurate than the AAM method, whatever the size of the sample. For example, for 100 simulated time series of length 200, with $H=0.65$, we obtain a standard deviation of the estimated $\theta$ very close to the average in the AAM case, whereas with ML the standard deviation is lower and the average closer to the actual value of $\theta$, as one can see in Table~\ref{tab:simulTheta}. We also note that reducing the sample size from 200 to 50 approximatively doubles the standard deviation of both the ML and AAM estimators.

\begin{table}[htb]
\centering
\begin{tabular}{|l|c|c|}
\hline
  Actual $\theta$ & AAM-estimated $\theta$ & ML-estimated $\theta$ \\
\hline
3 & 13.6  \textit{(8.3)} & 4.2 \textit{(5.7)} \\
10 & 26.8  \textit{(24.2)} & 10.8 \textit{(4.6)} \\
30 & 25.4  \textit{(19.6)} & 30.0 \textit{(11.5)} \\
50 & 45.4  \textit{(45.3)} & 52.0 \textit{(20.7)} \\
\hline
\end{tabular}
\begin{minipage}{0.7\textwidth}\caption{Average estimated $\theta$ parameter (and standard deviation) on 100 simulated time series of length 200 with the parameters $H=0.65$ and various $\theta$.}
\label{tab:simulTheta}
\end{minipage}
\end{table}

The results regarding the estimation of $H$ are still slightly in favour of the ML method as soon as the sample size is small, but the difference between the two methods is more restricted in particular for greater values of $H$, as one can see in Table~\ref{tab:simulH}. The real benefit of using the AAM method appears for greater sample sizes, for which the ML method is almost unusable due to unreasonably long execution time, as reported in Table~\ref{tab:execTime}. For the AAM approach, we get, in a limited time, results which have roughly the same accuracy as estimates obtained for shorter time series, as exposed in Table~\ref{tab:simulH2}.

\begin{table}[htb]
\centering
\begin{tabular}{|l|c|c|}
\hline
  Actual $H$ & AAM-estimated $H$ & ML-estimated $H$ \\
\hline
0.35 & 0.474  \textit{(0.240)} & 0.338 \textit{(0.102)} \\
0.50 & 0.556  \textit{(0.188)} & 0.510 \textit{(0.108)} \\
0.70 & 0.693  \textit{(0.174)} & 0.724 \textit{(0.126)} \\
0.80 & 0.774  \textit{(0.152)} & 0.838 \textit{(0.134)}  \\
\hline
\end{tabular}
\begin{minipage}{0.7\textwidth}\caption{Average estimated $H$ parameter (and standard deviation) on 100 simulated time series of length 50 with the parameters $\theta=30$ and various $H$.}
\label{tab:simulH}
\end{minipage}
\end{table}

\begin{table}[htb]
\centering
\begin{tabular}{|l|c|c|}
\hline
  $N$ & AAM execution time & ML execution time \\
\hline
1000 & 0.320 & 163.30 \\
500 & 0.377 & 20.24 \\
300 & 0.268 & 5.020 \\
200 & 0.216 & 1.704 \\
150 & 0.159 & 0.879 \\
100 & 0.119 & 0.336 \\
75 & 0.092 & 0.178 \\
50 & 0.084 & 0.075 \\
25 & 0.045 & 0.033  \\
\hline
\end{tabular}
\begin{minipage}{0.7\textwidth}\caption{Average execution time in seconds for the estimation of $H$ and $\theta$ with the AAM and ML methods for various lengths $N$ of time series.}
\label{tab:execTime}
\end{minipage}
\end{table}

\begin{table}[htb]
\centering
\begin{tabular}{|l|c|}
\hline
  Actual $H$ & AAM-estimated $H$  \\
\hline
0.40 & 0.444 \textit{(0.213)} \\
0.55 & 0.550  \textit{(0.262)} \\
0.65 & 0.635  \textit{(0.223)} \\
0.75 & 0.739  \textit{(0.205)} \\
0.80 & 0.786  \textit{(0.202)} \\
\hline
\end{tabular}
\begin{minipage}{0.7\textwidth}\caption{Average estimated $H$ parameter (and standard deviation) on 100 simulated time series of length 1,000 with the parameters $\theta=30$ and various $H$.}
\label{tab:simulH2}
\end{minipage}
\end{table}

The execution time displayed for the AAM method in Table~\ref{tab:execTime} in fact depends on some methodological choices, in particular the regression design. For instance, simulations show that the execution time depends linearly on the number $n$ of moments considered in the regression. This is consistent with the intuition. In Table~\ref{tab:execTime}, we have chosen $n=15$. However, beyond the execution time, the accuracy of the AAM method, both in terms of bias and variance, does not seem to be affected by $n$, according to simulations, at least for the values we have tested, between $n=15$ and $n=150$. The choice of the kernel and refinements on the choice of the bandwidth do neither seem to affect the accuracy.

\subsection{Estimation method and model specification}

We now show how the AAM framework makes it possible to determine whether the model is well specified or not. We known that the ML approach fails in this task, since it simply outputs parameters maximizing a likelihood, even when the model is not relevant. We can improve the ML approach by considering a set of various (fractal and stationary) models and ratios of likelihood or information criteria~\cite{BL}. This will in fact only help finding the best model in a predefined family of models but not clearly state the good specification of the selected model.

\begin{figure}[htbp]
	\centering
		\includegraphics[width=0.7\textwidth]{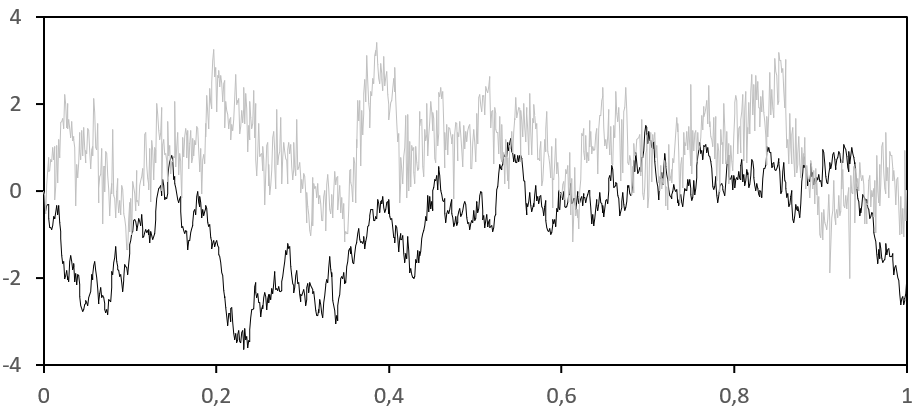} \\
		\includegraphics[width=0.5\textwidth]{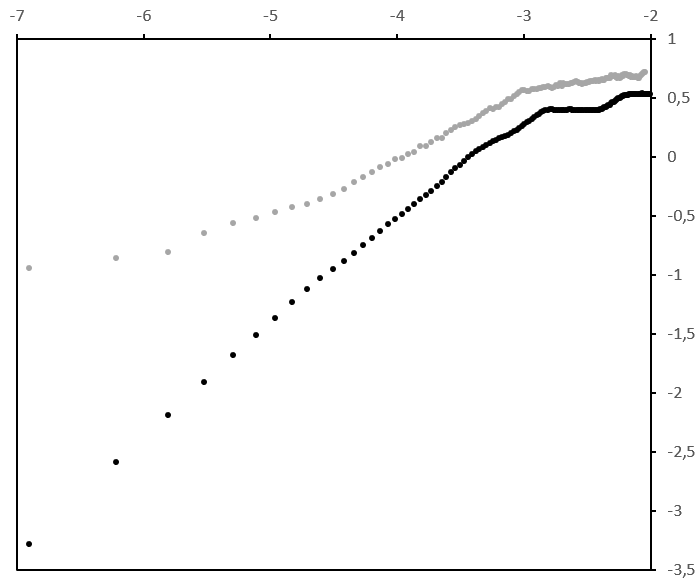}
\begin{minipage}{0.7\textwidth}\caption{Top: Simulated times series of a delampertized fBm (black) with the parameters $H=0.5$ and $\theta=30$, and a delampertized fBm of same parameters with an additive Gaussian noise of standard deviation $40\%$ (grey). Bottom: Corresponding log-log plot.}
	\label{fig:TimeSeries}
\end{minipage}
\end{figure}

We work with simulated data generated by two distinct models:
\begin{itemize}
\item \textbf{Time series 1} follows the specification mentioned above: it is a delampertized fBm of Hurst exponent $H=0.5$, variance 1, and time change parameter $\theta=30$. For this particular choice of Hurst exponent, the delampertized fBm is equivalent to an Ornstein-Uhlenbeck process.
\item \textbf{Time series 2} is the sum of two processes: a delampertized fBm with the same parameters as Time series 1, and a Gaussian noise process of standard deviation $40\%$.
\end{itemize}
We display the two trajectories of these processes in Figure~\ref{fig:TimeSeries}, as well as their corresponding log-log plot. For Time series 1, we observe a flattening of the curve at higher scales, which corresponds to the stationary feature. For Time series 2, in addition to the same flattening at high scales, we observe a second one at small scales. This is a well-known effect of additive noise on the log-log plot~\cite{LSG}.

We then use the ML approach to estimate the parameters. The likelihood used follows equation~\ref{eq:LL} and corresponds to a delampertized fBm. We thus expect this method to provide satisfying results for Time series 1 but not for Time series 2, since another specification is required. Table~\ref{tab:simulBruit} confirms this intuition with poor results for Time series 2, for instance a Hurst exponent equal to 0.267 very far from the true and unobserved value (0.5). However, when estimating the parameters with the ML method, we have no clue that the model and therefore the estimation are flawed. The estimated parameters are simply the best, provided that the model is well specified.

\begin{figure}[htbp]
	\centering
		\includegraphics[width=0.45\textwidth]{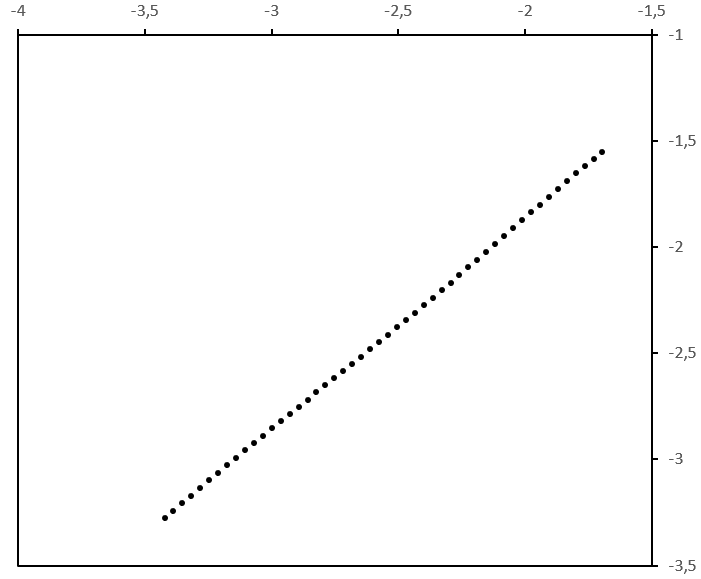}
		\includegraphics[width=0.45\textwidth]{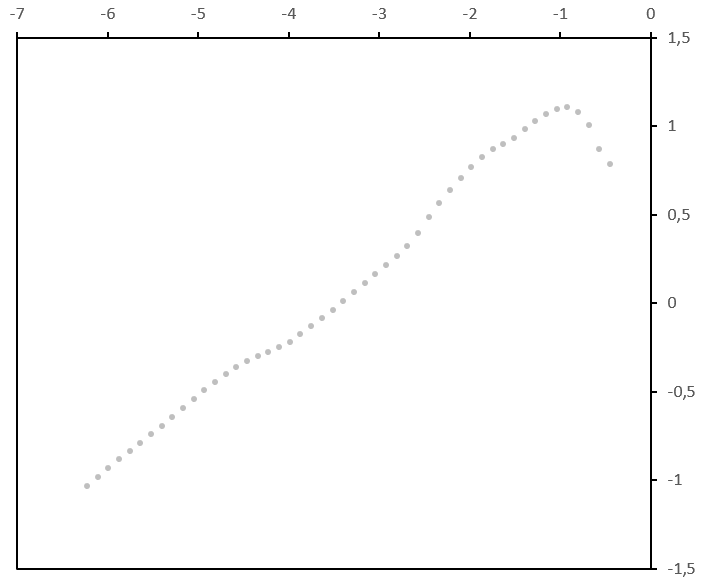}
\begin{minipage}{0.7\textwidth}\caption{Log-log plot of the Lamperti-transformation of Time series 1 (black) and Time series 2 (grey), using the parameters $H'$ and $\theta'$ estimated with the ML method with the assumption of a delampertized fBm.}
	\label{fig:LampLogplot}
\end{minipage}
\end{figure}

In order to check that the model is well specified, one can use the tools introduced for the AAM approach. More precisely, using the parameters $H'$ and $\theta'$ estimated with the ML method, one can Lamperti-transform the time series and consider the log-log plot of the transformed series, as displayed in Figure~\ref{fig:LampLogplot}. If the model is well specified, the log-log plot of the transformed series should be the one of an fBm, that is linear. For Time series 1, the log-log plot is indeed linear. Visually, we thus cannot assess that the delampertized fBm is not well specified for this time series. Regarding the transformed version of Time series 2, the aspect of the log-log plot is not really the one expected if the delampertized fBm was an appropriate choice of model, in particular because the curve stops growing for higher scales.

More quantitatively than visually, according to Theorem~\ref{th:LogPlotLin}, in case of good specification, the log-log plot of the transformed series should be linear with a slope equal to $2\times H'$. In this case, we thus should have $\hat H_{H',\theta'}$ close to $H'$ and $\hat{\alpha}_{H',\theta'}$ close to 1, where the metric $\hat H_{H',\theta'}$ is the half slope and $\hat{\alpha}_{H',\theta'}$ the linearity indicator introduced in Section~\ref{sec:AAM}. We see in Table~\ref{tab:simulBruit} that the linearity indicators of the two log-log plots are roughly equally close to 1. More interestingly, the quantity $\hat H_{H',\theta'}-H'$ is close to 0 (exactly 0.002) for Time series 1 and 36 times higher (0.072) for Time series 2. This undoubtedly confirms that a delampertized fBm is not an appropriate model for Time series 2.

\begin{table}[htb]
\centering
\begin{tabular}{|l|c|c|}
\hline
  & Time series 1 & Time series 2  \\
\hline
 $H'$: ML-estimated $H$ & 0.500 & 0.267 \\
 $\theta'$: ML-estimated $\theta$ & 32.1 & 1.97 \\
 $\hat H_{H',\theta'}$ & 0.498 & 0.195 \\
 $\hat{\alpha}_{H',\theta'}$ & 1.029 & 0.966 \\ 
\hline
\end{tabular}
\begin{minipage}{0.7\textwidth}\caption{Quantities of the ML and AAM approaches for both simulated time series.}
\label{tab:simulBruit}
\end{minipage}
\end{table}

\section{Financial application}

\begin{figure}[htbp]
	\centering
		\includegraphics[width=0.6\textwidth]{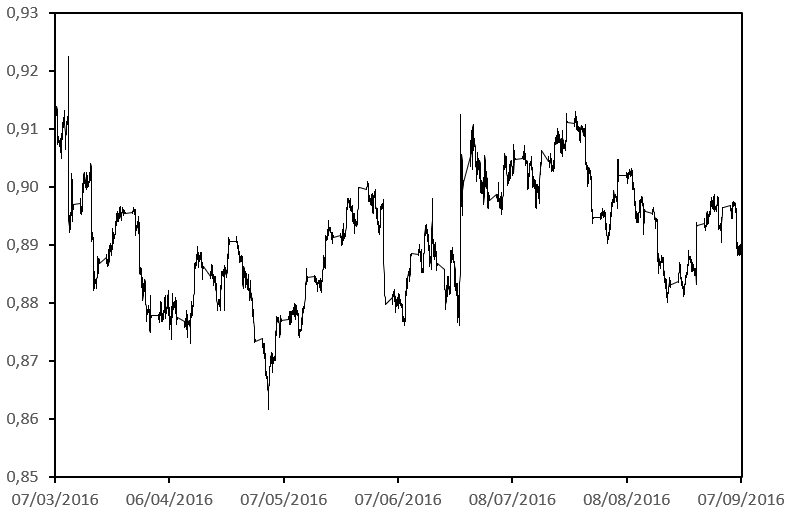}
\begin{minipage}{0.7\textwidth}\caption{Time series of log USD/EUR rate, sampled at a 15-minute time step.}
	\label{fig:FXseries}
\end{minipage}
\end{figure}

Following previous works on FX rates and Hurst exponents~\cite{DC,Garcin2017,GarcinLamperti}, we now focus on the time series of the logarithm of USD/EUR rate, sampled at a 15-minute time step, during the period starting the 7th March 2016 and finishing the 7th September 2016, as in Figure~\ref{fig:FXseries}. We split the sample in two parts of $6,460$ observations each. The first one is used for estimating models and the second half to test the predictive power of the models.

The analysis of the traditional log-log plot, between the scale of 15 minutes to the scale of 3 weeks, reveals a progressive decrease of the slope for higher scales, suggesting a stationarity of the time series. The perceived Hurst exponent for this range of scales is lower than $1/2$, more precisely it is $0.47$.

\begin{figure}[htbp]
	\centering
		\includegraphics[width=0.45\textwidth]{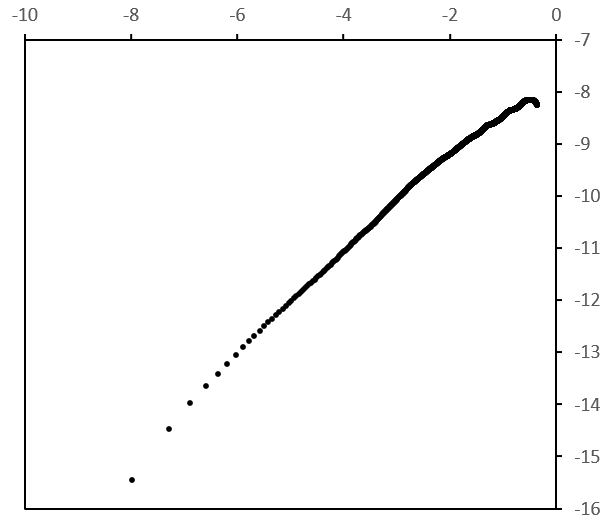}
		\includegraphics[width=0.45\textwidth]{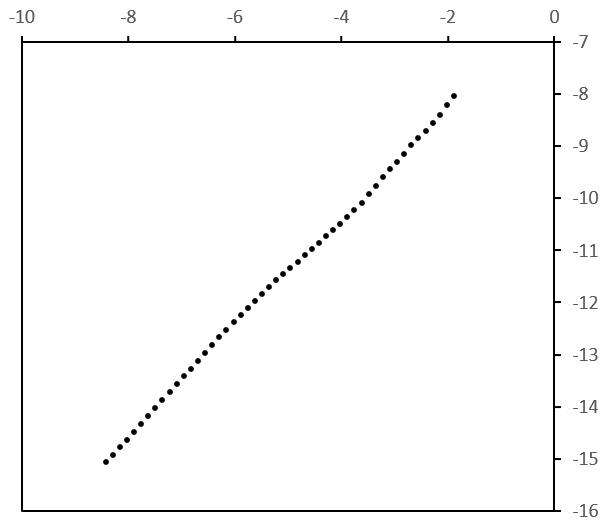}
\begin{minipage}{0.7\textwidth}\caption{Log-log plot of the series of log USD/EUR (left) and of its Lamperti transform (right), using the parameters $H'$ and $\theta'$ estimated with the AAM method with the assumption of a delampertized fBm.}
	\label{fig:FXLogplot}
\end{minipage}
\end{figure}

Applying the AAM estimator of the delampertized fBm leads to a linear log-log plot, with an $R^2$ of $0.999$ versus $0.982$ for the standard log-log plot. The Hurst exponent used for defining the Lamperti transform is found to be the same as half the slope of the log-log plot of the transformed series. This indicates the proper estimation of the model and of the transformation of a stationary process into a self-similar one. It is worth noting that the estimated underlying Hurst exponent is 0.53, that is above $1/2$ whereas the perceived Hurst exponent is below $1/2$. In other words, when one applies standard estimators of Hurst exponents on this series of FX rates, one concludes that the series is anti-persistent. But a finer analysis of the log-log plot, and the estimation of the Lamperti transform of the series, shows that this anti-persistence does not hold at every scale. It is indeed only the consequence of stationarity and price increments at small scales are in fact persistent, that is positively correlated. The two log-log plots, for the raw series and for its Lamperti transform, are displayed in Figure~\ref{fig:FXLogplot}. We note that the size of the dataset does not make it possible to use the ML approach, so that we have only used the AAM method.

In order to confirm the relevance of the model and of the estimator, we use the estimated Hurst exponents, the underlying and perceived ones, to build two competing forecasting procedures. Indeed, the value of the Hurst exponent is traditionally related to forecasting methods based on the covariance of the increments of an fBm or of any other selfsimilar model~\cite{NP,Mitra,ADG,GarcinForecast}. At each time $t$ in the test set, we forecast the variation of the FX rate between $t$ and $t+\tau$, where $\tau$ is the finest time scale of our sample, that is 15 minutes. We base our forecast on the observed  variation of the FX rate between $t-\tau$ and $t$. In other words, for a Hurst exponent above (respectively below) $1/2$, we anticipate same (resp. opposite) signs for the past and the future price increments. In particular, using the underlying Hurst exponent, we predict properly the sign of $56\%$ percent of the price increments in the test set. This ratio decreases to $44\%$ if we use instead the perceived Hurst exponent. With a confidence of more than $99.99\%$, a binomial test shows that these hit ratios are significantly different from $50\%$. 

\section{Conclusion}

In this paper, we have introduced two estimation methods for a delampertized fBm: the ML method and the AAM method. This second approach is an adaptation of the absolute-moment estimation method to the case of a stationary process. We have exposed the rationale of this method as well as a pseudo-code, and we have compared it to an ML method. The conclusion of this work highlights that the ML approach is in general more accurate than the AAM. But it also stresses the relevance of the AAM method compared to the ML in some cases, for the following reasons: it makes it possible to confirm visually that the model is well specified, thanks to the log-log plot, and it is computationally more performing. In practice, we recommend combining the two approaches: using the rapidity of the AAM approach to provide a first guess of the optimal parameters, to be used then for initializing the ML method, finally checking the good specification of the model with the log-log plot of a transformation of the data and with the statistics introduced in the AAM approach.


\bibliographystyle{plain}
\bibliography{bibliolamperti}

\begin{thebibliography}{10}

\bibitem{ADG}
A.~Ammy-Driss and M.~Garcin.
\newblock Efficiency of the financial markets during the {COVID}-19 crisis:
  time-varying parameters of fractional stable dynamics.
\newblock {\em arXiv preprint}, 2020.

\bibitem{ABKZ}
A.~Andresen, F.E. Benth, S.~Koekebakker, and V.~Zakamulin.
\newblock The {CARMA} interest rate model.
\newblock {\em International journal of theoretical and applied finance},
  17(2):1450008, 2014.

\bibitem{BLP}
M.~Bennedsen, A.~Lunde, and M.S. Pakkanen.
\newblock Decoupling the short-and long-term behavior of stochastic volatility.
\newblock {\em Journal of financial econometrics}, page nbaa049, 2021.

\bibitem{Bianchi}
S.~Bianchi.
\newblock Pathwise identification of the memory function of multifractional
  {B}rownian motion with application to finance.
\newblock {\em International journal of theoretical and applied finance},
  8(2):255--281, 2005.

\bibitem{BL}
S.~Bianchi and Q.~Li.
\newblock A new estimator of the self-similarity exponent through the empirical
  likelihood ratio test.
\newblock {\em Journal of statistical computation and simulation},
  90(11):1982--2001, 2020.

\bibitem{Box}
M.J. Box.
\newblock A new method of constrained optimization and a comparison with other
  methods.
\newblock {\em Computer journal}, 8(1):42--52, 1965.

\bibitem{BSZ}
D.C. Brody, J.~Syroka, and M.~Zervos.
\newblock Dynamical pricing of weather derivatives.
\newblock {\em Quantitative finance}, 2(3):189--198, 2002.

\bibitem{BI}
A.~Brouste and S.M. Iacus.
\newblock Parameter estimation for the discretely observed fractional
  {O}rnstein-{U}hlenbeck process and the {Y}uima {R} package.
\newblock {\em Computational statistics}, 28(4):1529--1547, 2013.

\bibitem{CKM}
P.~Cheridito, H.~Kawaguchi, and M.~Maejima.
\newblock Fractional {O}rnstein-{U}hlenbeck processes.
\newblock {\em Electronic journal of probability}, 8(3):1--14, 2003.

\bibitem{CTY}
Y.W. Cheung, D.C. Tam, and M.S. Yiu.
\newblock Does the {C}hinese interest rate follow the {US} interest rate?
\newblock {\em International journal of finance \& economics}, 13(1):53--67,
  2008.

\bibitem{Chevillard}
L.~Chevillard.
\newblock Regularized fractional {O}rnstein-{U}hlenbeck processes and their
  relevance to the modeling of fluid turbulence.
\newblock {\em Physical review {E}}, 96(3):033111, 2017.

\bibitem{CV}
A.~Chronopoulou and F.G. Viens.
\newblock Estimation and pricing under long-memory stochastic volatility.
\newblock {\em Annals of finance}, 8(2-3):379--403, 2012.

\bibitem{Coeur2001}
J.-F. Coeurjolly.
\newblock Estimating the parameters of a fractional {B}rownian motion by
  discrete variations of its sample paths.
\newblock {\em Statistical inference for stochastic processes}, 4(2):199--227,
  2001.

\bibitem{Coeur2005}
J.-F. Coeurjolly.
\newblock Identification of multifractional {B}rownian motion.
\newblock {\em Bernoulli}, 11(6):987--1008, 2005.

\bibitem{DC}
J.F. Diaz and J.H. Chen.
\newblock Testing for long-memory and chaos in the returns of currency
  exchange-traded notes ({ETN}s).
\newblock {\em Journal of applied finance and banking}, 7(4):15--37, 2017.

\bibitem{FBA}
P.~Flandrin, P.~Borgnat, and P.-O. Amblard.
\newblock From stationarity to self-similarity, and back: {V}ariations on the
  {L}amperti transformation.
\newblock In {\em Processes with long-range correlations}, pages 88--117.
  Springer, Berlin-Heidelberg, 2003.

\bibitem{Garcin2017}
M.~Garcin.
\newblock Estimation of time-dependent {H}urst exponents with variational
  smoothing and application to forecasting foreign exchange rates.
\newblock {\em Physica {A}: statistical mechanics and its applications},
  483:462--479, 2017.

\bibitem{GarcinLamperti}
M.~Garcin.
\newblock Hurst exponents and delampertized fractional {B}rownian motions.
\newblock {\em International journal of theoretical and applied finance},
  22(5):1--26, 2019.

\bibitem{GarcinMulti}
M.~Garcin.
\newblock Fractal analysis of the multifractality of foreign exchange rates.
\newblock {\em Mathematical methods in economics and finance}, 13-14(1):49--73,
  2020.

\bibitem{GarcinForecast}
M.~Garcin.
\newblock Forecasting with fractional {B}rownian motion: a financial
  perspective.
\newblock {\em arXiv preprint}, 2021.

\bibitem{GG}
M.~Garcin and M.~Grasselli.
\newblock Long versus short time scales: the rough dilemma and beyond.
\newblock {\em to appear in {D}ecisions in economics and finance}, 2020.

\bibitem{GJR}
J.~Gatheral, T.~Jaisson, and M.~Rosenbaum.
\newblock Volatility is rough.
\newblock {\em Quantitative finance}, 18(6):933--949, 2018.

\bibitem{HSL}
S.~Hidot, C.~Saint-{J}ean, and J.-Y. Lafaye.
\newblock Etude exp\'erimentale de l'influence d'un \'echantillonnage
  irr\'egulier dans l'estimation du param\`etre de {H}urst.
\newblock {\em Journal de la soci\'et\'e fran\c{c}aise de statistique},
  149(1):81--95, 2008.

\bibitem{HN}
Y.~Hu and D.~Nualart.
\newblock Parameter estimation for fractional {O}rnstein-{U}hlenbeck processes.
\newblock {\em Statistics \& probability letters}, 80(11-12):1030--1038, 2010.

\bibitem{Lamperti}
J.~Lamperti.
\newblock Semi-stable stochastic processes.
\newblock {\em Transactions of the {A}merican mathematical society},
  104(1):62--78, 1962.

\bibitem{LSG}
Y.~Lanoisel\'ee, G.~Sikora, A.~Grzesiek, D.S. Grebenkov, and
  A.~Wy\l{}oma\'nska.
\newblock Optimal parameters for anomalous-diffusion-exponent estimation from
  noisy data.
\newblock {\em Physical review {E}}, 98(6):062139, 2018.

\bibitem{LF}
F.~Le~Floc'h.
\newblock Issues of {N}elder-{M}ead simplex optimisation with constraints.
\newblock {\em S{SRN} e{L}ibrary}, page 2097904, 2012.

\bibitem{MvN}
B.~Mandelbrot and J.~van {N}ess.
\newblock Fractional {B}rownian motions, fractional noises and applications.
\newblock {\em S{IAM} review}, 10(4):422--437, 1968.

\bibitem{MW}
J.~Mielniczuk and P.~Wojdy\l{}\l{}o.
\newblock Estimation of {H}urst exponent revisited.
\newblock {\em Computational statistics \& data analysis}, 51(9):4510--4525,
  2007.

\bibitem{Mitra}
S.K. Mitra.
\newblock Is {H}urst exponent value useful in forecasting financial time
  series?
\newblock {\em Asian social science}, 8(8):111--120, 2012.

\bibitem{MM}
M.~Mossberg and E.~Mossberg.
\newblock A note on parameter estimation in {L}amperti transformed fractional
  {O}rnstein-{U}hlenbeck processes.
\newblock {\em I{FAC} proceedings volumes}, 45(16):1067--1072, 2012.

\bibitem{NS}
P.K. Narayan and S.S. Sharma.
\newblock Does data frequency matter for the impact of forward premium on spot
  exchange rate?
\newblock {\em International review of financial analysis}, 39:45--53, 2015.

\bibitem{NM}
J.A. Nelder and R.~Mead.
\newblock A simplex method for function minimization.
\newblock {\em Computer journal}, 7(4):308--313, 1965.

\bibitem{NP}
C.J. Nuzman and H.V. Poor.
\newblock Linear estimation of self-similar processes via {L}amperti's
  transformation.
\newblock {\em Journal of applied probability}, 37(2):429--452, 2000.

\bibitem{PLV}
R.F. Peltier and J.~L\'evy~{V}\'ehel.
\newblock A new method for estimating the parameter of fractional {B}rownian
  motion.
\newblock {\em Technical report 2396, {INRIA}}, 1994.

\bibitem{RL}
W.~Robbertse and F.~Lombard.
\newblock On maximum likelihood estimation of the long-memory parameter in
  fractional {G}aussian noise.
\newblock {\em Journal of statistical computation and simulation},
  84(4):902--915, 2014.

\bibitem{Shi}
L.~Shi.
\newblock Does security transaction volume-price behavior resemble a
  probability wave?
\newblock {\em Physica {A}: statistical mechanics and its applications},
  366:419--436, 2006.

\bibitem{ST}
T.~Simos and M.~Tsionas.
\newblock Bayesian inference of the fractional {O}rnstein-{U}hlenbeck process
  under a flow sampling scheme.
\newblock {\em Computational statistics}, 33(4):1687--1713, 2018.

\bibitem{Viitasaari}
L.~Viitasaari.
\newblock Representation of stationary and stationary increment processes via
  {L}angevin equation and self-similar processes.
\newblock {\em Statistics \& probability letters}, 115:45--53, 2016.

\bibitem{SGKMBP}
M.~\v{S}apina, M.~Garcin, K.~Kramari\'c, K.~Milas, D.~Brdari\'c, and
  M.~Piri\'c.
\newblock The {H}urst exponent of heart rate variability in neonatal stress,
  based on a mean-reverting fractional {L}\'evy stable motion.
\newblock {\em Fluctuation and noise letters}, 19(3):2050026, 2020.

\end{thebibliography}

\appendix

\section{Proof of Theorem~\ref{th:MomentIncr}}\label{sec:ProofMomentIncr}

\begin{proof}
The results concerning $X$ and $Y$ are standard properties~\cite{FBA,GarcinLamperti}. Let us now focus on the process $Z$.
\begin{itemize}
\item First, we note that, whatever $t> 0$, 
\begin{equation}\label{eq:TransTrans}
Z_t=\left(\mathcal L_{H',\theta'}\mathcal L_{H,\theta}^{-1}X\right)_t=t^{h}X_{t^{\theta/\theta'}},
\end{equation}
where we introduced $h=H'-\frac{\theta}{\theta'}H$. This process is extended by $Z_0=0$. We obtain this extension trivially if $h>0$, and using equation~\eqref{eq:TransTrans} else: $\E[Z_t^2]=t^{2h}t^{2H\theta/\theta'}=t^{2H'}\overset{t\rightarrow 0}{\longrightarrow} 0$.
\item The process $Z$ has Gaussian increments of mean 0 and of variance:
\begin{equation}\label{eq:VarIncr}
\begin{array}{ccl}
\E\left[(Z_{t+\tau}-Z_{t})^2\right] & = & (t+\tau)^{2h}\E\left[X_{(t+\tau)^{\theta/\theta'}}^2\right] + t^{2h}\E\left[X_{t^{\theta/\theta'}}^2\right] - 2(t+\tau)^{h}t^h\E\left[X_{(t+\tau)^{\theta/\theta'}}X_{t^{\theta/\theta'}}\right] \\
 & = & \sigma^2\left[ (t+\tau)^{2H'} + t^{2H'} - (t+\tau)^{h}t^h\left((t+\tau)^{2H\theta/\theta'}+t^{2H\theta/\theta'}-\left|(t+\tau)^{\theta/\theta'}-t^{\theta/\theta'}\right|^{2H}\right) \right],
\end{array}
\end{equation}
where we used successively equation~\eqref{eq:TransTrans} and equation~\eqref{eq:CovFBM}. 
\item This variance of the increments varies with $t$. We prove this property by considering the asymptotic behaviour of the variance, for $\tau/t\rightarrow 0$. Indeed, we have, using equation~\eqref{eq:VarIncr} together with asymptotic expansions and the fact that $H<1$:
\begin{equation}\label{eq:VarLim0}
\begin{array}{ccl}
\E\left[(Z_{t+\tau}-Z_{t})^2\right] & = & \sigma^2 \left[(t+\tau)^{2H'}\left(1-\left(\frac{t}{t+\tau}\right)^h\right)+t^{2H'}\left(1-\left(\frac{t+\tau}{t}\right)^h\right) + (t+\tau)^{h}t^h\left|(t+\tau)^{\theta/\theta'}-t^{\theta/\theta'}\right|^{2H}\right] \\
& \overset{\tau/t\rightarrow 0}{\sim} & \sigma^2 \left[h\frac{\tau}{t}\left((t+\tau)^{2H'}-t^{2H'}\right) + t^{2h} t^{2H\theta/\theta'}\left|1+\frac{\theta}{\theta'}\frac{\tau}{t}-1\right|^{2H}\right] \\
& \overset{\tau/t\rightarrow 0}{\sim} & \sigma^2 \left[2hH't^{2H'}\left(\frac{\tau}{t}\right)^2 + t^{2H'}\left(\frac{\theta}{\theta'}\frac{\tau}{t}\right)^{2H}\right] \\
& \overset{\tau/t\rightarrow 0}{\sim} & \sigma^2t^{2H'}\left(\frac{\theta}{\theta'}\frac{\tau}{t}\right)^{2H}.
\end{array}
\end{equation}
This asymptote depends on $t$, as soon as $H\neq H'$, so increments are not stationary in this case. If $H=H'$, we cannot conclude that increments are stationary and we must go further in our stationarity analysis. We thus consider the variance of increments with a given $\tau>0$, for $t=0$: 
\begin{equation}\label{eq:VarLim0bis}
\E\left[(Z_{0+\tau}-Z_{0})^2\right] =\E\left[(Z_{\tau})^2\right]=\sigma^2\tau^{2H'},
\end{equation}
and for $t\rightarrow+\infty$, that is for $\tau/t\rightarrow 0$, using equation~\eqref{eq:VarLim0} with $H=H'$:
\begin{equation}\label{eq:VarLim0ter}
\underset{t\rightarrow +\infty}{\lim}\ \E\left[(Z_{t+\tau}-Z_{t})^2\right] = \sigma^2\left(\frac{\theta}{\theta'}\tau\right)^{2H'}.
\end{equation}

The two variances, when $H=H'$, are equal only if $\theta=\theta'$. We thus have the necessary condition: stationarity of increments is only possible when both $H=H'$ and $\theta=\theta'$. 

The sufficient condition is straightforward. Indeed, when $H=H'$ and $\theta=\theta'$, then $h=0$ and equation~\eqref{eq:VarIncr} gives exactly, $\forall t>0$, $\E\left[(Z_{t+\tau}-Z_{t})^2\right]=\sigma^2\tau^{2H}$, which does not depend on $t$.

\item Increments of $Z$ are Gaussian variables whose variance is provided by equation~\eqref{eq:VarIncr}. Then, using the definition of $M_{k,N,t_a,t_b}(Z)$ and the fact that the absolute moment of a standard Gaussian variable $G$ is $\E\left[|G|^k\right]= 2^{k/2}\Gamma\left(\frac{k+1}{2}\right) / \Gamma\left(\frac{1}{2}\right)$, we get:
$$\E\left[M_{k,N,t_a,t_b}(Z)\right]= \frac{A(\sigma,k)}{N}\sum_{i=1}^{N}{\left[ t_{i+1}^{2H'} + t_i^{2H'} - t_{i+1}^{h}t_i^h\left(t_{i+1}^{2H\theta/\theta'}+t_i^{2H\theta/\theta'}-\left[t_{i+1}^{\theta/\theta'}-t_i^{\theta/\theta'}\right]^{2H}\right) \right]^{k/2}},$$
which is the result displayed in Theorem~\ref{th:MomentIncr}.

\item For the asymptotic value of $\E\left[M_{k,N,t_a,t_b}(Z)\right]$, we use equation~\eqref{eq:VarLim0}, with the notation $\tau=(t_b-t_a)/N$:
\begin{equation}\label{eq:Th1Asym1}
\begin{array}{ccl}
\E\left[M_{k,N,t_a,t_b}(Z)\right] & \overset{\tau\rightarrow 0}{\sim} & \frac{A(\sigma,k)}{N}\sum_{i=1}^{N}{\left[ \left(\frac{\theta}{\theta'}\tau\right)^{2H}t_i^{2(H'-H)} \right]^{k/2}} \\
 & \overset{N\rightarrow +\infty}{\sim} & \frac{A(\sigma,k)}{N}\left(\frac{\theta}{\theta'}(t_b-t_a)/N\right)^{kH}\sum_{i=1}^{N}{t_i^{k(H'-H)}}.
\end{array}
\end{equation}
We recognize the Riemann sum:
\begin{equation}\label{eq:Th1Asym2}
\begin{array}{ccl}
\frac{1}{N}\sum_{i=1}^{N}{t_i^{k(H'-H)}} & \overset{N\rightarrow +\infty}{\sim} & \frac{1}{t_b-t_a}\int_{t_a}^{t_b}{t^{k(H'-H)}dt} \\
 & \overset{N\rightarrow +\infty}{\sim} & \frac{1}{t_b-t_a}\frac{t_b^{k(H'-H)+1}-t_a^{k(H'-H)+1}}{k(H'-H)+1}.
\end{array}
\end{equation}
Equations~\eqref{eq:Th1Asym1} and~\eqref{eq:Th1Asym2} together lead to the asymptotic expression in Theorem~\ref{th:MomentIncr}.

\end{itemize}
\end{proof}

\section{Proof of Theorem~\ref{th:LogPlotLin}}\label{sec:ProofLogPlotLin}

\begin{proof}
We prove first the sufficient condition. If $(H',\theta')=(H,\theta)$, then $Z=X$, $f_{\mathcal S}(H',\theta')=0$ for all $\mathcal S$, and the statement is trivial.

Regarding the necessary condition, if we assume that $(H',\theta')$ reaches the theoretical minimum of $f_{\mathcal S}$, which is zero according to the particular case $(H',\theta')=(H,\theta)$, then $\alpha_{H',\theta'}=1$ and $\hat H_{H',\theta'}=H'$. As we have the same $\alpha_{H',\theta'}$ and $\hat H_{H',\theta'}$ whatever $\mathcal S$, we conclude that the plot $\ln(\tau)\mapsto \ln\left(M_{H',\theta',\tau}\right)$ is affine. Regarding the slope, after equation~\eqref{eq:VarLim0}, it is $2H$ for small scales, but also for all scales because the log-log plot is affine, so that $\hat H_{H',\theta'}=H$. The value of $H'$ which minimizes $f_{\mathcal S}$ is thus $H$. If we follow the notation of Theorem~\ref{th:MomentIncr}, with $H'=H$, we should have $\ln\left(M_{2,N,t_a,t_b}(Z)\right)\overset{N\rightarrow+\infty}{\sim}\ln\left(M_{2,1,t_a,t_b}(Z)\right)-2H\ln(N)$, that is, following equations~\eqref{eq:th1NonAsympt} and~\eqref{eq:th1Asympt}, 
$$2H\left(\ln\left(\frac{\theta}{\theta'}\right)+\ln(t_b-t_a)\right)=\ln\left(t_b^{2H}+t_a^{2H}-(t_bt_a)^{H(1-\theta/\theta')}\left(t_b^{2H\theta/\theta'}+t_a^{2H\theta/\theta'}-(t_b^{\theta/\theta'}-t_a^{\theta/\theta'})^{2H}\right)\right),$$
which is true only if $\theta'=\theta$. Therefore $(H',\theta')=(H,\theta)$.
\end{proof}

\end{document}